\documentclass[12pt]{article}

\usepackage{epsfig,amsmath,latexsym}
\usepackage{amsfonts,amssymb,graphicx}
\usepackage{amscd, amssymb, amsthm, amsmath,eucal, verbatim}

\newtheorem{theorem}{Theorem}[section]

\newtheorem{lemma}[theorem]{Lemma}

\newcommand{\ds}{\displaystyle}
\newcommand{\bx}{{\bf x}}

\newcommand{\RR}{{\mathbb{R}}}
\newcommand{\sG}{{\mathcal{G}}}
\newcommand{\sH}{{\mathcal{H}}}

\begin{document}

\title{Area Inequalities for Embedded Disks Spanning Unknotted Curves}

\author {Joel Hass$^1$, Jeffrey C. Lagarias$^1$ and William P.
Thurston$^2$}

\footnotetext[1]{The first two authors carried out some of this work
while visiting the Institute for Advanced Study.
The first author was partially supported by NSF grant DMS-0072348,
and by a grant to the Institute by AMIAS.}

\footnotetext[2]{Partially supported by NSF grant DMS-9704286.}
\maketitle
\date{}

\begin{abstract}
We show that a smooth unknotted curve in $\RR^3$ satisfies
an isoperimetric inequality that bounds the area of an embedded
disk spanning the curve in terms of two parameters: the length
$L$ of the curve and the thickness $r$ (maximal radius of an
embedded tubular neighborhood) of the curve.
For fixed length, the expression giving the
upper bound on the area grows exponentially in $1/r^2$.
In the direction of lower bounds, we
give a sequence of length one curves with $r \to 0$ for which
the area of any spanning disk is bounded from
below by a function that grows exponentially with $1/r$.
In particular, given any constant $A$,
there is a smooth, unknotted length one curve for which the
area of a smallest embedded spanning disk is greater than $A$.
\end{abstract}


\noindent

%
%
%
%

\section{Introduction}
\setcounter{equation}{0}

The classical isoperimetric inequality gives an upper  bound for the
area of the disk bounded by a curve in $\RR^2$.
For a simple closed curve of length $L$
in the plane that encloses area $A$, it states that
\begin{equation}~\label{isop}
4 \pi  A \le L^2,
\end{equation}
with equality if and only if the curve is a circle.
There are many extensions of the isoperimetric
theorem to higher dimensions; some of these
are discussed in Almgren~\cite{Alm86}, Burago and Zalgaller~\cite{BZ88},
and Osserman \cite{Os78}.

There are two natural extensions of the isoperimetric inequality
to a simple closed $C^2$-curve $\gamma$ in $\RR^3$ of length $L$.

\noindent
(1) There exists an immersed smooth disk
in $\RR^3$, with $\gamma$ as  boundary, having area $A$
with  $ 4\pi A \le  L^2.$ This disk may self-intersect and may
intersect the curve $\gamma$, as indeed it must if the curve is
knotted. If the curve is not a circle, then there exists
such an immersed disk for which strict inequality holds.

\noindent
(2) There exists an  orientable smooth
surface embedded in $\RR^3$, with no restriction on its genus,
with $\gamma$ as boundary, having area $A$,
with $4\pi A \le  L^2.$
If the curve is not a circle, then there exists
such a surface for which strict inequality holds.

Result (1) traces back to a result of Andre Weil in 1926
concerning isoperimetric inequalities  for surfaces of
non-positive Gaussian curvature.
See also Beckenbach and Rado~\cite{BR33}.
Result (2) traces back to a 1930
argument in Blaschke \cite[p. 247]{Bla30},
see Osserman \cite[p. 1202]{Os78}.
For completeness we indicate  proofs of these well-known 
results in an Appendix (\S\ref{classical}).
We derive (2) directly from (1).
The hypotheses in (1) and (2) can be weakened to require only
that $\gamma$ be rectifiable; however we assume $C^2$-curves
for ease of comparison with later results.

This paper treats the
situation when the two conditions above are imposed simultaneously.
That is, we suppose the curve $\gamma$ is unknotted, and therefore
is the boundary of an embedded disk, and ask for upper
bounds on the area  of some embedded spanning disk. 
Note that an embedded spanning disk attaining the  minimal area 
under these conditions need not exist, 
because the embeddedness condition may not be preserved under limits.
However there is still an infimum for the area of all embedded spanning disks,
which we refer to as the minimal area of a spanning disk.

We will show that under these conditions there is
no upper bound on the minimal area
given in terms of the length $L$ of $\gamma$ alone.
To obtain an upper bound on the area of some spanning
disk we must therefore impose more geometrical
control on  $\gamma$ than just bounding its length.
For this we  consider an
additional geometric quantity, the thickness
of $\gamma$.

The {\em thickness} $r(\gamma)$ of a $C^2$-curve
$\gamma$ in $\RR^3$ is
the supremum of $r$ such that a radius $r$
normal tubular neighborhood of the curve is embedded.
This notion can be used to define a knot invariant.
The supremum 
of  the ratio $r(\gamma)/L(\gamma)$
over all $C^2$-curves $\gamma$ having a given knot type 
defines a knot invariant
called the {\em thickness} of a knot, studied in
Buck and Simon \cite{BS99} and
Litherland et al. \cite{LSDR99}, where it is
related to various knot energies. 
For a $C^2$-closed curve $\gamma$
the  quantity
$L(\gamma)/r(\gamma)$ is called the {\em ropelength}
of $\gamma$, and the infimum of the  ropelength
over all $C^2$-representatives of a knot type,
called the minimal ropelength, 
is the inverse of the thickness.
Various properties of ropelength, including the
existence of a minimizer in the class $C^{1,1}$, have been
investigated recently  \cite{CKS02}.

Our first two results concern area bounds in terms of
the length $L$ and thickness $r$ of an
unknotted  closed curve $\gamma$.
In \S\ref{upper.bounds2} we show that for an unknotted curve
there is an upper bound for the area of a spanning disk
depending on these quantities. \\

%
%

\begin{theorem}~\label{upperbound0}
For any embedded closed unknotted $C^2$-curve $\gamma$ in $\RR^3$ having
length $L$ and thickness $r$,
there exists a smooth embedded disk of area $A$  having
$\gamma$ as boundary with
$$
A \le (C_0)^{(L/r)^2} L^2,
$$
where $C_0 > 1$ is a constant independent of $\gamma$, $L$ and $r$.
\end{theorem}

In the other direction, 
we show that an upper bound of this
type must be at least exponential in $L/r$, and that
there is no upper bound on the area depending solely on $L$. \\

%
%

\begin{theorem} \label{noupperbound1}
There is a constant  $C_1>1$
and a sequence of unknotted, $C^{\infty}$-curves
$\{\Delta_n\}_{n=1,\dots,\infty}$ embedded in $\RR^3$,
each having  length $L_n=1$ and thickness $r_n > 1/n$,
with the property that for large enough $n$,
any embedded disk spanning $\Delta_n$ has area
$$
A \ge (C_1)^{n} \ge (C_1)^{L_n/r_n} L_n^2  .
$$
\end{theorem}

We include the constant $L_n = 1$ for comparision to Theorem~\ref{upperbound0}.
These exponential lower bounds differ from the
exponential upper bound in Theorem~\ref{noupperbound1}
in having an exponent linear in $L_n/r_n$ rather
than quadratic. Such examples necessarily have
$1/n \le r_n \le c/ \sqrt{n}$ for a constant $c$ as $n \to \infty$,
by comparing the lower bound $(C_1)^n$ with
the upper bound in  Theorem~\ref{upperbound0}.

We next show that a bound on the length
together with 
 an upper bound on the curvature of
an unknotted curve is also
not sufficient to give an upper bound for
the area of some spanning disk. Without loss of 
generality we may
normalize by setting the length of the curve to one. \\

%
%

\begin{theorem}~\label{noupperbound2}
There is an infinite sequence
$\{\Gamma_n\}_{n=1,\dots,\infty}$ of unknotted $C^\infty$-curves
of length $L_n=1$ embedded in $\RR^3$, each with (pointwise) curvature
everywhere bounded above by a constant $K_0$, independent of $n$,
such that any embedded disk spanning $\Gamma_n$ has area
$$
A \ge n.
$$
\end{theorem}

The proofs of Theorems~\ref{noupperbound1} and \ref{noupperbound2}
are both given in \S\ref{lower.bounds}.
They are proved via a detour
into piecewise-linear topology, and an investigation of
similar questions for piecewise-linear curves, i.e. embedded polygons
in $\RR^3$. For piecewise linear curves, one can consider
combinatorial analogues of isoperimetric inequalities, in
which the ``length'' is replaced by the number of edges $n$ in
a polygon, and the ``area'' is the number of triangles in
an embedded triangulated surface having the polygon as boundary, as in 
\cite{HassLagarias:02}.
We establish  the following PL analogue of Theorem~\ref{upperbound0}
above. \\

%
%

\begin{theorem}~\label{PLupperbound00}
There is a constant $C_2 >1$ such that
given any  polygonal  unknotted curve $P$ in $\RR^3$ containing
$n$ line segments, there is a triangulated PL disk embedded in
$\RR^3$ having $P$ as boundary that contains at most $(C_2)^{n^2}$
triangles.
\end{theorem}

This result is established  as part of Theorem~\ref{PLupperbound}
in \S2. To deduce Theorem~\ref{upperbound0} from it, we show that
we can obtain such a triangulated surface in which the area of
each triangle is bounded above by a constant.

Hass, Snoeyink and Thurston
\cite{HST} established the following
PL analogue of Theorems~\ref{noupperbound1} and~\ref{noupperbound2}. \\

%
%

\begin{theorem}~\label{PLowerbound00}
There is a constant $C_3 >1$ and an infinite sequence of
unknotted polygonal curves
$\{K_n\}_{n=1,...,\infty}$ having between $n$ and $23n$ segments,
such that any  triangulated PL disk embedded in
$\RR^3$ having $K_n$ as boundary contains at least $(C_3)^{n}$ triangles.
\end{theorem}

This result is not quite sufficient to  deduce
Theorem~\ref{noupperbound1}.
It is  necessary to strengthen the construction in \cite{HST}, and
we carry this out  in Theorem~\ref{properties} in \S4.
Using this strengthened  construction, in \S5 we show that
a large number of triangles occurring in the triangulated PL disk have
their
areas bounded below by a fixed positive constant, and that such
the resulting lower bound carries
over to smooth embedded disks with the same
or nearby boundaries.

In the direction of PL combinatorial bounds,
Hass and Lagarias \cite{HassLagarias:02} showed
that there is a (qualitative) combinatorial
analogue of the isoperimetric inequality (\ref{isop}) above
in $\RR^3$.

The plan of the paper is as follows.
In  \S\ref{upper.bounds} we derive
the PL upper bound  Theorem \ref{PLupperbound00}. From it
we obtain upper bounds on the area of
a piecewise-smooth embedded disk spanning a
polygonal unknot $P$ with $n$ edges. In \S\ref{upper.bounds2} we
deduce Theorem~\ref{upperbound0} from these results.
In \S\ref{PL} we review some  results of \cite{HST},
and construct a family of unknotted polygonal curves 
$\{K_n\}_{n=1,...,\infty}$ that do
not have  PL spanning disks having few triangles,
and have additional useful properties.
The families of curves $\{\Delta_n\}_{n=1,...,\infty}$ and 
\{$\Gamma_n\}_{n=1,...,\infty}$ used to prove
Theorems~\ref{noupperbound1} and \ref{noupperbound2}
are constructed by a smoothing process from these curves.
In  \S\ref{lower.bounds} we relate the PL complexity to area,
and obtain lower bounds for the area of any
spanning disk of $\Gamma_n$, resp.  $\Delta_n$, that
increase to infinity with $n$, and deduce
Theorems~\ref{noupperbound1} and \ref{noupperbound2}.

We note that the  proofs of the lower bound estimates
Theorems~\ref{noupperbound1} and \ref{noupperbound2}
in \S\ref{PL} and \S\ref{lower.bounds}
can be read independently of  the proofs of the
upper bound results Theorem~\ref{upperbound0}
and \ref{PLupperbound00} in \S\ref{upper.bounds} and
\S\ref{upper.bounds2}. 

%
%
%
%

\section{Upper bound:  Minimal complexity of a PL
embedded spanning disk}~\label{upper.bounds}

Let $P$ be a polygonal, unknotted curve of length $L$ in
$\RR^3$ having $n$ edges.
In this section we deduce upper bounds on the number
of triangles in some embedded
PL spanning disk for $P$, and
for the minimal area of such a disk.

%
%

\begin{theorem} \label{PLupperbound}
There exists a constant  $C_4 > 1$ such
that whenever $P$ is a closed, unknotted, polygonal curve
in $\RR^3$ having $n$ edges and length $L$, the following hold:

\noindent
(1) There exists an embedded PL disk spanning $P$
whose number of triangles is at most
$\displaystyle \frac{1}{32} (C_4)^{n^2}$.

\noindent
(2) Such an embedded PL disk
can be chosen to lie inside a ball of radius $4L$, and
its area then satisfies
$$
A \leq (C_4)^{n^2}L^2.
$$
\end{theorem}

Theorem~\ref{PLupperbound} is proved by a modification of
ideas used in Hass and Lagarias \cite{HassLagarias:01}.
That paper obtains upper bounds on the number of elementary
moves required to move from an unknotted polygonal curve
in the one-skeleton of a triangulated 3-manifold having $t$ tetrahedra,
to a single triangle, by a sequence of moves which
deform the curve across three
surfaces, an annulus, a torus and a disk. Here we need to produce a
spanning disk. The construction of \cite{HassLagarias:01}
can be modified to produce such a disk, with the bound on the
number of triangles in the disk having the same functional form
${C}^t$ as the bounds obtained in \cite[Theorem 1.2 ]{HassLagarias:01},
but with a different constant $C > 1$. A second necessary change
is that we must construct a triangulation of a 3-ball containing the
polygon $P$ in its one-skeleton, and we do not have the freedom
to move this curve; our construction obtains a quadratic upper bound
for the number of tetrahedra in the constructed 3-ball triangulation.

%
%
 
\begin{theorem} \label{unknot.PLmanifold}
There exists a constant $C_5 > 1$ such that
whenever $M$ is a triangulated PL $3$-manifold
with boundary, containing $t$ tetrahedra, and
$K$ is a closed polygonal curve contained in the $1$-skeleton
of the interior of $M$ that is unknotted in
$M$, then there is a PL triangulated disk in $M$ that
has a subdivision of $K$ as boundary, and that
contains at most $(C_5)^t$ triangles.
\end{theorem}

{\bf Proof:}
Theorem~1.2 of \cite{HassLagarias:01} states that there is a constant
$c_2$, independent of $M$,
such that if $K$ is an unknot embedded in the 1-skeleton of
interior$(M)$, then $K$ can be isotoped to a single triangle
using at most $2^{c_2 t}$ elementary moves.
The proof of this result deforms the original knot $K$ to
a triangle with three sequences of elementary moves,
made inside three triangulated
surfaces.
A peripheral torus is constructed by taking two
barycentric subdivisions and taking the torus that is
the boundary of a regular neighborhood $R_K$ of $K$.
An embedded disk is
used to construct elementary moves taking
an essential curve $K_2$ on $\partial R_K$ to a single triangle in $M$.
On $\partial R_K$ there is a curve $K_1$ with prescribed number of segments
that is a longitude, and is isotopic to $K$ in $R_K$ by an
isotopy with an explicitly computed bound on the number of elementary
moves.
An embedded annulus is used to go from $K$ to the longitude
$K_1$. The 2-torus $\partial R_K$ is then used to define a
sequence of elementary moves from $K_1$ to the
boundary of a normal disk $K_2$. These moves define a homotopy between
two embedded curves on $\partial R_K$.

The disk and the annulus used to construct the elementary moves are
each embedded in $M$.
However the elementary moves
that connect $K_1$ to $K_2$ in $\partial R_K$
may go back and forth over the torus $\partial R_K$
many times, and the annulus that is swept out on $\partial R_K$
may not be embedded.
This problem is rectified by a modification of the construction given in
\cite{HassLagarias:01} that perturbs this swept out
annulus to an embedding.

In Theorem~1.2 of \cite{HassLagarias:01} it is shown that:

\noindent
(1) A normal disk with a longitude $K_2$ as boundary contains at most
$2^{4608t+6}$ triangles.

\noindent
(2) There is an embedded $PL$ annulus $S$
in the solid torus $R_K$ whose two boundary components are
another longitude $K_1$ and $K$,
and which consists of at most $2^{1858t}$ triangles.

\noindent
(3) The number of elementary moves needed to deform $K_1$ to $K_2$ on
the torus $\partial R_K$ is at most $ 2^{10^7t}$.

Let $n$ denote the number of
elementary moves taking $K_1$ to $K_2$, so that $n \le  2^{10^7t}$ and
let $T$  denote the PL surface $\partial R_K$.
We retriangulate $M$ so as to
insert a product $ T \times [0,1]$ in place of
the triangulated torus $T$, with
$T \times 0$ and $T \times 1$ triangulated identically to $T$.
The interior of $ T \times [0,1]$ is triangulated as follows:
The torus $T \times t_i$ is triangulated identically to $T \times 0$, $i=
1, \dots, n$.
We triangulate the product $T \times [t_i, t_{i+1}]$,  $i= 1, \dots, n-1$
in two steps. Begin by dividing this region into
triangular prisms, one for each
triangle on $T \times t_i$. Then divide
each prism into 14 tetrahedra, by first dividing
each rectangular face of the prism into four triangles with a common
vertex at
the center of the rectangle, and then coning the 14 triangles on the
boundary
of the prism to a new vertex at the center of the prism.
We now construct an embedded annulus in $ T \times [0,1]$
whose boundary on $T \times 0$ is $K_2$ and whose
boundary on $T \times 1$ is $K_1$.
The intersection of the annulus
with $T \times (t_i, t_{i+1})$ consists
of the image of $K_1$ in $T$ after $i-1$ elementary moves
$\times (t_i,t_{i+1})$.
On $T \times t_i$ the annulus consists of
the image of $K_1$ in $T$ after $i-1$ elementary moves plus the single
triangle
swept out by the $i$th elementary move if the $i$th move is of type 2.
If the  $i$th move is of type 1 then it consists of
{image of $K_1$ in $T$ after $i-1$ elementary moves}
plus or minus a single vertex corresponding to the $i$th move.
The image of $K_1$ on $T \times t_i$ after $i$ elementary moves
for any $i$ consists of less than
$ 2^{10^7t}2^{1858t} =  2^{10^7t+1858t}$ edges.
Each vertical annular piece lying in $T \times (t_i, t_{i+1})$
can be triangulated with
four triangles for each vertical rectangle above an edge
in $T \times t_i$, and each
horizontal piece lying in a torus $T \times t_i$ contains at most one
triangle.
There are at most $ 2^{10^7t}$ such annuli in $ T \times [0,1]$,
one for each elementary move.
So the total number of triangles required to triangulate the
annulus in $T \times [0,1]$ is at most $ (2^{10^7t+1858t + 2})(2^{10^7t})$.
We obtain a total of at most $(2^{10^7t})(2^{10^7t+1858t + 2}) \le
2^{10^8t-1}$
triangles in the annulus connecting $K_2$ on $T \times 0$
to $K_1$ on $T \times 1$.

Adding to this the contributions of the disk and
annulus in (1) and (2) above results in
a spanning disk containing less than $2^{10^8t}$ triangles.
We take $C_5 =  2^{10^8}$.
$~~~\qed$

To apply this result, we construct a triangulated convex polyhedron $B$
in $\RR^3$ which contains the polygonal curve $P$ in its
one-skeleton; $B$ will be the PL $3$-manifold in the result above.

%
%

\begin{lemma} \label{triangulation}
Let $P$ be an unknotted polygonal curve in $\RR^3$ having
$n$ vertices and length $L$.
Then there exists a triangulated convex polyhedron  $B$ (with
interior vertices) in $\RR^3$, which contains $P$ strictly in its
interior, in its $1$-skeleton, and which has the
following properties.

(1) The manifold $B$ is contained in a sphere of radius $4L$ in $\RR^3$.

(2) The number of tetrahedra in $B$ is less than $290n^2+ 290n +116$.

\end{lemma}
{\bf Proof:}
The proof is similar in spirit to the construction of
a triangulated $3$-manifold given in \cite[Theorem 8.1]{HassLagarias:01}.
However the argument in \cite{HassLagarias:01} must
be modified, because here the polygonal
curve $P$ is fixed, and its vertices must be in
the triangulation, whereas in \cite{HassLagarias:01} only the projection
is specified, and there is freedom
to move around the vertices to simplify the triangulation.

After rescaling $P$ by the homothety $ \displaystyle \bx \to
\bx /L$,
we may assume its length is $1$.
In this case $P$  lies in the interior of
a ball of radius $1/2$ in $\RR^3$, which can
without loss of generality be assumed to be centered at the origin.
By rotating the coordinate system if necessary, we can suppose
that the projection of $P$ in the $z$-direction onto
the $xy$-plane is a regular projection, and then,
as in Lemma 7.1 of \cite{HLP2}, this projection
has less than $n^2$ crossings.
The graph associated to the projection has straight-line
edges, and its vertices consist of the projections of
vertices of the polygon $\gamma$ together with all crossing points.

As a preliminary step, we construct an augmented version  of this
projection graph in the $z = 0$ plane, which will be
a triangulated planar graph that we call $\sG$.
Let $S'$ denote the convex hull of the vertices of the projection
graph, and add
extra edges as necessary to fill out the boundary
of the convex hull. Note that all crossing vertices in
the projection lie strictly in the interior of the convex
hull. Add new vertices to the graph
by adding  on  each edge of this graph two new interior
``special'' vertices, spaced one-third of the way along the
edge from its two endpoints. The resulting graph has at
most $5n^2 + 5n$ vertices. This graph is planar with convex
faces, and we triangulate it as follows. We first enclose each
crossing vertex in four triangles, by adding the four edges
connecting pairs of interior ``special'' vertices nearest
to it.  We then triangulate all remaining
(convex) polygonal faces of the graph in an arbitrary fashion,
adding no new vertices. As a final step we add a border to this graph
consisting of three points arranged in an equilateral
triangle in the $z=0$ plane at distance $3$ from the origin,
together with the three edges of this triangle. Note that
each edge gets no closer to the origin than $\displaystyle
3\sqrt{3}/2$,
and so remains outside the convex hull $S'$ of the projection graph,
which is contained in the ball of radius $1/2$ around
the origin. Finally we triangulate the
outer region thus added by connecting each vertex
on the boundary of $S'$ to appropriate vertices of the
equilateral triangle.

The resulting triangulated graph $\sG$
has at most $5n^2 + 5n + 3$ vertices, and so there  are exactly
$10n^2 + 10n + 4$ triangles in its triangulation. Its convex hull
is contained in the circle of radius $3$ centered at
the origin.

To begin construction of the triangulated polyhedron,
we start with two vertically translated
copies of this augmented graph,  in the planes
$z = 1$ and $z = -1$. The convex hull of this set is
a triangular prism, which will be the polyhedron $B$.
We proceed to triangulate $B$, i.e. to dissect it into
tetrahedra sharing common faces, such that $\gamma$ is in the
one-skeleton. We will use the triangulation of the
upper and lower faces on $z= 2$ and $z = -2$ arising
from the triangulation of the augmented graph.

We proceed to construct a simplicial complex  $\sH$ inside
$B$ which will include $P$ in its 1-skeleton, and
which vertically projects to $\sG$ on the plane
$z = 0$; we will include $\sH$ in the triangulation and not include
$\sG$.
The simplicial complex $\sH$ is mostly two-dimensional,
consisting of triangles lying over each triangle of $\sG$ not
including a crossing vertex, but thickened
to a three-dimensional
tetrahedron over each triangle of $\sG$ that does include
a crossing vertex. To begin, we add
the polygon $P$, and add as new vertices along
each edge of $P$ the points that vertically project to
crossing vertices of the graph; there are then two
vertices on $P$ that project to each crossing vertex
of the graph. To these upper and lower crossing vertices of $P$
add a vertical edge connecting them, which is called
the  vertical edge corresponding to the crossing.
Then add extra vertices on this configuration,
two along each edge,
that project to the ``special'' vertices of the graph
$\sG$. We add the three equilateral vertices of
$\sG$ lying in the plane $z = 0$, plus the  edges
connecting them to the vertices in $\sH$ that project
to the corresponding vertices in $\sG$ that they connect to.
Each vertex in $\sG$ that is not a crossing vertex has
a unique ``lifted'' vertex in $\sH$ that vertically projects to it,
while crossing vertices have two ``lifted'' vertices, an upper
and a lower one.
Each triangle in the triangulation $\sG$ not containing
a crossing vertex then lifts uniquely to a triangle in $\sH$
having the corresponding ``lifted'' vertices.
Corresponding to each crossing vertex of $\sG$ we create
four tetrahedra in the simplicial complex $\sH$, which share the
vertical edge corresponding to the crossing, and completely
surround it. Each of the four tetrahedra has as vertices
the upper and lower crossing vertex of $P$ plus
two ``special'' vertices whose projection gives a
triangle in the triangulation $\sG$. We have now specified
the simplicial complex $\sH$. It has the property that
each top-dimensional simplex in it projects vertically to
a unique triangle in $\sG$ and each such triangle is covered once.
Furthermore  $\sH$ has an ``upper surface'' and a
``lower surface'' obtained by extracting from each tetrahedron
its triangular upper and lower face having the same vertical
projection.

We now add vertical edges connecting each vertex in the
``upper surface'' of $\sH$ to the corresponding
vertex of the graph $\sG$ in the plane $z=2$. We can match
triangular faces in this copy of $\sG$ with the triangular
face in the ``upper surface'' and their convex hull is
a convex triangular prism, with boundary edges the vertical
edges we just added. Each such triangular prism may be
subdivided into fourteen  tetrahedra, by
first subdividing each of its three
quadrilateral vertical faces into four
triangles by adding a vertex at its
barycenter, and then coning all the
triangular faces to the barycenter of
the triangular prism. Triangulations
of adjacent prisms are automatically compatible (adjoining triangular
faces match); this triangulates the ``upper'' portion
of $B.$ The ``lower'' portion is cut up into compatible
tetrahedra in an exactly similar fashion, adding vertical
edges. In this triangulation of $B$, we have $28$ tetrahedra
corresponding to  each triangle in $\sG$ that does not
include any crossing vertex, and $29$ tetrahedra for those
faces that do. Thus this triangulation uses
at most $290n^2 + 290n + 116$ tetrahedra.

We have obtained a
linear triangulation of a convex region $B$
contained in the ball of radius $4$ such that
$P$ lies in the one-skeleton of $B$, strictly in
the interior of $B$. Upon rescaling by the homothety $\bx \to L\bx$
we obtain the desired triangulation, with $B$ inside a sphere of
radius $4L$. $~~~\qed$ \\

\noindent{\bf Remark.} We do not know whether the
bound $O(n^2)$ in the triangulation above
is the optimal order of magnitude.
Avis and ElGindy \cite{AE87} give examples of sets of
$n$ points in $\RR^3$, for arbitrarily large $n$,
such that any triangulation of the convex hull
including these points as vertices required at least
$\Omega(n^2)$ tetrahedra. However, they  also showed that
for $n$ points in general position (no four in a plane)
there exists a triangulation of the convex hull using
$O(n)$ tetrahedra. We do not know whether such an improvement
to $O(n)$ is possible in the case above for general position vertices,
because we impose the stronger requirement that $P$
be in the $1$-skeleton of the triangulation.
See~\cite{Ch84} for a possible approach. \\

{\bf Proof of Theorem ~\ref{PLupperbound}:}
(1) Let $B$ be the triangulated ball
constructed in Lemma~\ref{triangulation},
that has at most $t = 290n^2 + 290n + 116$ tetrahedra.
By hypothesis $P$ is
unknotted in $\RR^3$, and it follows that it is
unknotted in the ball $B$, since an unknotting disk
can be moved into $B$ by an isotopy. By Theorem~\ref{unknot.PLmanifold}
there is an embedded PL spanning disk in $B$ having $P$
as boundary and containing at most $(C_5)^t$ triangles. This gives
an upper bound of $\displaystyle \frac{1}{32} (C_4)^{n^2}$ triangles, with
$\displaystyle (C_4)= 32(C_5)^{696}$, on noting that $n \ge 3$.
This yields (1).

(2) The embedded PL spanning disk constructed in
Lemma~\ref{triangulation} is contained in a ball $B$ of radius
$4L$. Each triangle in it therefore has area at most $32L^2$.
Thus the total area is bounded by
$\displaystyle \frac{1}{32}(C_4)^{n^2} (32L^2),$
which yields (2). $~~~\qed$ \\

\noindent
{\bf Remark:} We have spanned the polygonal curve $P$ with
a PL disk. Such disks can be approximated on their
interior by smooth disks of
no greater area, so Theorem~\ref{PLupperbound} implies the same
upper bound for the area of a smooth spanning disk for $P$.

%
%
%
%

\section{Upper bound: Minimal area  of
a smooth embedded spanning disk} ~\label{upper.bounds2}
In this section we consider closed unknotted $C^2$-curves $\gamma$
embedded in $\RR^3$.
We obtain upper  bounds for the minimal area of an embedded smooth disk
having the curve as boundary as a function of the length $L$
and  the thickness
$r$ of the curve. In \S\ref{lower.bounds} we will
show that there is no
upper bound for the minimal area of an embedded disk
having the curve as boundary that is purely a function of
the length and the maximal curvature of such a curve.

Let $\gamma$ be a rectifiable closed curve in $\RR^3$ of length
$L$, and let $\{\gamma(s):~ 0 \le s \le L\}$
be a parameterization of the curve
by arc length, starting from a fixed point $\gamma(0)$ on the curve.
Recall that an {\em  $r$-neighborhood} $T_{r}$ for
a curve $\gamma$ consists of all points in $\RR^3$ within
Euclidean distance $r$ of some  point on the curve.
An $r$-neighborhood $T_r$  is {\em tubular}
if it is foliated by disjoint two-dimensional closed
Euclidean disks $D_s(r)$ of radius $r$ isometrically
embedded in $\RR^3$, whose center is the point $\gamma(s),$
for $0 \le s < L$, and each disk is normal to the curve
$\gamma$ at the point $\gamma(s)$.
If so, we call the disks $D_s(r)$ {\em fibers} of the
$r$-tubular neighborhood.
Every closed $C^2$-curve embedded in $\RR^3$
has an $r$-tubular neighborhood for some positive $r$,
see  \cite[Theorem 5.1, pp. 109-110]{Hi76}.

%
%
\begin{theorem} \label{unknot.smoothupper}
Suppose that $\gamma$ is an unknotted $C^2$-curve
in $\RR^3$ of length $L$ that has an $r$-tubular neighborhood.
Then there is an embedded smooth disk having boundary $\gamma$ whose
area $A$ satisfies
$$
A < (C_6)^{(L/r)^2} L^2,
$$
with $C_6$ is a constant independent of $\gamma$.
\end{theorem}

The idea of the proof is to approximate $\gamma$ by an inscribed polygonal
closed curve $P'$ which remains entirely
inside a thinner normal tubular neighborhood  and
has $n \le 32L/r + 1$ segments.
By Theorem~\ref{PLupperbound} there exists a
PL spanning disk for $P'$ having area at most $(C_4)^{n^2} L^2$.
This is at most $C^{(L/ r)^2} L^2$, for a different
constant $C$.
We deform this PL spanning disk to a spanning disk for $\gamma$
using a homeomorphism $\tau$ of $\RR^3$ that moves $P'$ to $\gamma$,
is the identity outside the
$r$-tubular neighborhood of $\gamma$, and is a piecewise
differentiable Lipschitz map
with Lipschitz constant everywhere bounded by 5.
It follows that the area of the resulting spanning disk for $\gamma$
is larger by at most a factor of $25$ over
that of the PL spanning disk for $P'$.  The Lipschitz map
we construct is the identity
outside a neighborhood of $\gamma$; this neighborhood 
is  a sort of necklace of beads around $\gamma$, in which each
bead is a right circular cone, with cone point on
the curve   $\gamma$,  rather then a round ball.
We first establish some preliminary results concerning
the deformations of a curve to  an inscribed polygon.

%
%

\begin{lemma}~\label{cone}
Let $\alpha: [0,ar] \to \RR^3 $ be a smooth, properly embedded unit speed
$C^2$-curve with curvature bounded above by $1/r$, with $0 < a \le 1$.
The tangent vectors $T(s)$ of $\alpha$ lie inside a 
right circular cone with cone point at $\alpha(0)$,
axis on the line in direction  $\alpha '(0)$, and with
cone angle $ \theta_a = 2\arcsin {a/2}$ subtended to its axis.
The entire curve $\alpha(s)$ lies within this closed cone.
\end{lemma}
{\bf Proof:}  Let $T(s)$ denote the unit tangent at $0 \le s \le ar$ and
$k(s)=||T'(s)||$ the magnitude of the curvature vector. Then
$||T'(s)|| =  k(s) \le 1/r$, and
$$
||T(s) - T(0)|| = || \int_0^{s} T'(t) dt || \le \int_0^{ar} ||T'(t)|| dt
\le (ar)(1/r) = a.
$$
It follows from the geometry of the isosceles triangle formed by
$T(0), T(s),$ and $T(s) - T(0)$ 
that all tangent vectors  to $\alpha(s) $ lie in a cone of angle
$ \theta_a = 2\arcsin {a/2}$ centered around the line through $T(0)$.
In particular,
$$
T(0) \cdot T(s) \ge \cos ( 2\arcsin {a/2}) > 0
$$
over this whole range, since $0 < a \le 1$.
The curve $\alpha$ itself also lies inside
the cone with axis in the direction  $T(0)$ and based at $\alpha(0)$.
For if the curve ever  left this cone then it would also leave a cone with a
slightly larger cone angle, and would exit the boundary of this slightly
wider cone  transversely.  The tangent
vector of the $\alpha$ at such a point would not lie
 within the original cone. $~~~\qed$

We next compare the distance of two points on a curve with specified tubular
neighborhood as
measured along the curve and as measured in $\RR^3$.

%
%
\begin{lemma}~\label{relative distance}
Let $0< a < 1/2$  and let $\alpha:[0,ar] \to \RR^3 $ be a smooth,
properly embedded unit speed
$C^2$-curve with curvature everywhere bounded by $1/r$. 
Let $e$ denote the line
segment in $\RR^3$ connecting $\alpha(0)$ and $\alpha(ar)$. Then
$\alpha$ is contained in the right circular cone with 
cone point $\alpha(0)$,  axis $e$, and that subtends the angle
$2 \theta_a = 4 \arcsin (a/2)$ to its axis. Moreover
$$
(\mbox{length}(e))^2  \ge \cos( 2 \theta_a)~ (\mbox{length}(\alpha))^2.
$$
\end{lemma}

{\bf Proof:}
Lemma~\ref{cone} 
shows that  two tangent vectors on $\alpha(s)$ differ by an angle
at most $2 \theta_a = 4 \arcsin (a/2)$, so that their inner 
product satisfies
$T(s_1) \cdot T(s_2) \ge \cos (2\theta_a) > 0$,
using $0 < a \le 1/2$ in the last inequality. 
Let $T$ denote the unit tangent vector pointing from $\alpha(0)$
to $\alpha(ar)$, so that 
$$
T = \frac{1}{\mbox{length}(e)}(\alpha(ar) - \alpha(0)) =
 \frac{1}{\mbox{length}(e)}   \int_{0}^{ar} T(s) ds.
$$
Then for any $0 \le s' \le ar$, the inner product 
\begin{eqnarray*}
T \cdot T(s') &= &  \frac{1}{\mbox{length}(e)}~
\int_{0}^{ar} T(s) \cdot T(s') ds  \\
& \ge & \frac{ar}{\mbox{length}(e)}~\cos (2 \theta_a) \\
& \ge &  \cos(2 \theta_a). 
\end{eqnarray*}
It follows that the entire curve lies in the right circular
cone with vertex at $\alpha(0)$, axis the line segment $e$
pointing in direction $T$,
subtending angle $2\theta_a = 4 \arcsin(a/2)$, and with its
endpoint $\alpha(ar)$ being the center of the
base of the cone.

Now since $\mbox{length}(\alpha) = ar$ we obtain
\begin{eqnarray*}
\mbox{length}(e) & = & ||\alpha(ar) - \alpha(0)|| = 
|\int_{0}^{ar} T(s) \cdot T(s') ds| \\
& = &   | \int_{0}^{ar} T(s') \cdot \left(    
\frac{1}{\mbox{length}(e)}\int_{0}^{ar} T(s) ds \right) ds' | \\
& \ge &\frac{1}{\mbox{length}(e)}\int_{0}^{ar}\int_{0}^{ar}
|T(s') \cdot T(s) | ds ds' \\
& \ge & \frac{ (\mbox{length}(\alpha))^2}{\mbox{length}(e)} \cos (2 \theta_{a}).
\end{eqnarray*}
This yields 
$$(\mbox{length}(e))^2 \ge  
\cos (2 \theta_{a})(\mbox{length}(\alpha))^2, 
$$ 
as asserted.
$~~~\qed$

In what follows we let $\gamma$ denote a 
$C^2$-curve having a normal tubular neighborhood of radius $r$.
The existence of the tubular neighborhood for $\gamma$ 
implies that the absolute
value of the curvature 
at each point of $\gamma$ is at most $1/r$.
(The points where the normal exponential map develops singularities
are called the focal set.  They occur at distance $r$ for a curve at 
a point  of (absolute) curvature $1/r$, see \cite[p. 232]{DC92}.)
Therefore Lemmas~\ref{cone} and  ~\ref{relative distance} are applicable.

We construct a particular  polygonal curve $P'$  approximating the 
$C^2$-curve $\gamma$ as follows.
Let  $P'$  in $\RR^3$ have vertices on the
curve $\gamma$ taking $\ds n  = \lfloor\frac{32L}{r}\rfloor  +1$ 
equally spaced points
$z_0, z_1, \dots z_{n-1}$ along $\gamma$ (measured by
arclength of $\gamma$) and inscribing line segments connecting
successive points. Let
$\gamma_j$ be the arc of $\gamma$ between $z_j$ and $z_{j+1}$ and 
let $e_j = [z_j, z_{j+1}]$ be the line segment of
$P'$ between these points in $\RR^3$, with last segment
$e_{n-1} = [z_{n-1}, z_0].$  
All arcs $\gamma_j$ are of equal length $L/n$ and we define the
quantity $a_j$ by $a_jr= \mbox{length}(\gamma_j) = L/n$, which implies
that $1/33 \le a_j \le 1/32.$ We will verify
below that $P'$ is embedded in $\RR^3$.

Let $C_j$ denote the right circular cone with cone point 
$z_j$, axis the line segment $e_j$, subtending the
angle $4 \theta_{a_j} = 8 \arcsin (a_j/2)$,
Note that $C_j$ has base diameter no larger than
$2 \sin(4\theta_{a_j}) (a_j r)$ 
since $\mbox{length}(e_j) \le a_jr$,
and the center of its base is $z_{j+1}$.

%
%

\begin{lemma}~\label{cones.disjoint}
Let $\gamma$ be a closed $C^2$-curve with tubular neighborhood
of radius $r$. In the construction above, the arc $\gamma_j$
lies inside the closed cone $C_j$.
Distinct cones $C_j$ and $C_k$ with $j<k$ are disjoint
unless $j=k+1$ or $j=n-1, k=0$, and in this case the two cones
intersect in the single point $z_{j+1}$, resp. $z=0$ if
$j=n-1$.
The polygonal curve $P'$ is embedded in $\RR^3$.
\end{lemma}

{\bf Proof:} Lemma~\ref{relative distance}
applies to the arc $\gamma_j$ 
to show that it is contained inside
the cone  $C_j$.
The cone $C_j$ subtends an
angle to its axis of 
$$
4\theta_{a_j} = 8 \arcsin a_j/2 \le 8 \arcsin (1/64) < 0.14,
$$
and the base of this cone has diameter at most
$$
2 \sin (4 \theta_{a_j})(a_j r) \le 2 (8 \arcsin (1/64)) r/32
< 0.01 r.
$$ 
If two cones $C_j$ and $C_k$ intersect, then they necessarily contain
points $v_j \in \gamma_j$ and $v_k \in \gamma_k$ separated by
distance at most $(0.02) r$. Let $e$ be the line segment
between $v_j$ and $v_k$.

We claim that since $\gamma$ has an $r$-tubular 
neighborhood, 
the shorter of the two arcs of $\gamma$ connecting $v_j$ and $v_k$ 
has length at most $(0.03)r$. 
To see this, consider the arc of $\gamma$ starting
from $v_j= \gamma(s_0)$ to a variable point $w= \gamma(s')$ 
where $s'$ is an arclength parametrization of $\gamma$ that
increases from $s = s_0$ as $\gamma(s')$ moves along the
curve in one direction. The distance $||v_j - \gamma(s')||$ is
initially increasing. Let $\gamma(s_2)$ be the first maximum point reached.
At that point the tangent vector $T(s_2)$ is orthogonal to the
line $[v_j, w]$, and part of this line lies in the
normal disk of radius $r$ around $w$ so
we  must have $||v_j - w|| \ge r$. Similarly allowing $s'$
to decrease from $s$ we find that $|| \gamma(s')- v_j||$
increases monotonically to the first point $s_1$ where  
$||v_j - \gamma(s_1)||$
reaches a maximum, and this maximum is at least $r$ by the same
reasoning. 

Now the remaining part of 
the curve falls outside the normal tubular neighborhood of radius $r$ determined by
the arc $\{ \gamma(s) : s_1 \le s \le s_2\}$.
This normal tubular neighborhood includes  all
points within distance $r/4$ of $v_j$, for 
if $v$ is a point with $||v - v_j|| \le r/4$
then the closest point $w= \gamma(s^{\ast})$ on
the arc $\{ \gamma(s) : s_1 \le s \le s_2\}$
to $v$ has $||w-v|| \le r/4$
from it, and attains a local minimum for the distance function from $v$.
It is at distance at least $r/2> r/4$
from each of the endpoints $\gamma(s_1), \gamma(s_2)$,
since otherwise 
$$||v_j - \gamma(s_1)|| \le ||v_j - w|| + ||w - v||
+ ||v -\gamma(s_1)|| < r/4 + r/4 + r/2 = r,
$$
a contradiction.
The local minimality implies that
$v-w$ is perpendicular to the tangent vector
$T(s^{\ast})$, so lies in the disk of radius
$r$ at $w= \gamma(s^{\ast})$ normal to $T(s^{\ast})$. 
We conclude that all
points on the curve $\gamma$
within distance $r/2$ of $v_j$ must lie in that
part of the curve where $||v_j - \gamma(s')||$ is monotone
increasing as $|s'-s_0|$ increases. The remark
above shows that this holds at least
for $|s' - s_0| \le r$.
Since $||v_j - v_k || < 0.02r$, the point $v_k$ must fall on this
part of the curve. Now we can apply
Lemma~~\ref{relative distance} to observe
that on taking $a= 0.03$
the segment $e= [\gamma(0) , \gamma( 0.03r)]$ already satisfies
$$
(\mbox{length}(e))^2 \ge \cos (2 \theta_a)(ar)^2
\ge \cos (4 \arcsin 0.015)(0.03 r)^2 > (0.02)^2 r^2.
$$
We conclude that the shorter of the two arcs of 
$\gamma$ connecting $v_j$ and $v_k$ 
has length at most $(0.03)r$, as asserted.

Since the length of each segment
$\gamma_j$ is at least $r/33 > 0.03 r$, it follows that
$C_j$ and $C_k$ are adjacent cones, i.e. $k = j+1$ or $j=n-1, k=0$.
Adjacent cones however
intersect in a single point, the center of the base of one and
the cone point of the other. 

Finally, since each segment $e_j$ of $P^{'}$ lies in  a different
cone, and overlaps are only possible at endpoints of the $e_j's$,
the closed polygonal curve $P'$ is embedded in $\RR^3$.
$~~~\qed$

The next two lemmas will be used
to construct a Lipschitz homeomorphism inside
the tubular neighborhood that carries $P'$ to $\alpha$.
We start with a two-dimensional map.
Given a two-dimensional Euclidean disk $D(d) =D(0, d)$ of
radius $d$ centered at $(0,0)$ and a point
$q= (x,y)$ in the disk, with $||q|| < d$,
let $\sigma_{q,d} : D(d) \to D(d)$ be 
the map that sends the origin to $q$ and that maps each
line segment
from  the origin to  a boundary point of $D(d)$ linearly to the
line segment from $q$ to the same boundary point. It is 
explicitly given by
$$
\sigma_{q,d}(w) = w + (1-  ||w||/d)q.
$$

%
%
%
\begin{lemma}~\label{sigma.lipschitz}
The map $\sigma_{q,d}: D(d) \to D(d)$ is a Lipschitz homeomorphism
with Lipschitz constant at most $ \displaystyle 1 + ||q||/d \le 2$.
Its inverse map $\sigma_{q,d}^{-1}$ is Lipschitz
with Lipschitz constant at most
$ \displaystyle  \frac{1}{1 - (||q||/d)}.$ If $q_1$ and $q_2$ lie in
the disk  $D(d)$, then for any point $w \in D(d)$
$$
||\sigma_{q_1,d}(w) - \sigma_{q_2,d}(w)|| \le ||q_1 - q_2||.
$$ 
\end{lemma}
{\bf Proof:}
The map $\sigma_{q,d}$ is clearly a homeomorphism,
which leaves the boundary of $D(d)$ fixed.  It remains to
estimate a Lipschitz bound.
For $w \in D(d)$ we have
$$
\sigma_{q,d}(w) = w + (1-||q||/d)q.
$$
We calculate that
\begin{eqnarray*}
|| \sigma_{q,d}(w_2) - \sigma_{q,d}(w_1)|| & \le & || w_1 -w_2|| + 
\frac{||q||}{d}(||w_1||-||w_2||) \\
& \le & ||w_1 -w_2||(1 + ||q||/d).
\end{eqnarray*}
Similarly $ \displaystyle || \sigma_{q,d}(w_2) - \sigma_{q,d}(w_1)|| \ge 
||w_1 -w_2||(1 - ||q||/d)$.
So  $\sigma_{q,d}$ has Lipschitz constant bounded above by 
$1+||q||/d $
and below by  $1 - ||q||/d$.
It follows that the inverse map $\sigma_{q,d}^{-1}$ has Lipschitz constant
bounded above by $\displaystyle  \frac{1}{1 - (||q||/d)}.$

Finally, we have
$$
||\sigma_{q_1,d}(w) - \sigma_{q_2,d}(w)|| =
|(1-  ||w||/d)|||q_1 - q_2|| \le ||q_1 - q_2||,
$$
as asserted. $~~~\qed$ \\

Next we give a parametrized generalization of 
the map of Lemma~\ref{sigma.lipschitz}
which is a three-dimensional map.
We define in $(x, y, t)$-space the right circular cone 
$$
C(b,h) := \{ (x,y,t) : x^2 + y^2 \le b^2t^2, ~~0 \le t \le h\}
$$
of  height $h$ and base of radius $bh$, centered on the $t$-axis,
and suppose that $0 < b \le 1.$
We consider a curve $\alpha(t)= \{( x(t), y(t), t ) : 0 \le t \le h\}$ 
parametrized by the $t$-variable 
which is  strictly monotone increasing in $t$, 
with
$\alpha(0) = (0,0,0) $ and $\alpha(h) = (0, 0, h)$.
and we assume that the whole curve 
is contained in  the  right circular cone $C(b,h)$.
Concerning the smoothness of this curve we assume only that 
 $t \to \alpha(t)$ is Lipschitz-continuous with
Lipschitz constant $C_L(\alpha)$, i.e. that 
$$
||\alpha(t_1) - \alpha(t_2)|| \le C_L(\alpha)|t_1 - t_2|
$$
holds for $0 \le t_1, t_2 \le h.$
 
Now, associated to this curve  $\alpha$, we
 define a homeomorphism  $\tau = \tau_{\alpha}$ of $\RR^3$ supported in 
the cone $C(b,h)$, which leaves its boundary fixed.
It takes   horizontal disks with fixed $z$-coordinate into themselves
moving each of them by the map
$\sigma_{\alpha(t),bt}$ of Lemma~\ref{sigma.lipschitz}, and moves
the point $(0, 0, t)$ on
the axis of the cone to $\alpha(t)= (x(t), y(t), t)$.
This homeomorphism is defined inside the cone $C(b, h)$ by
$$
\displaystyle \tau (x,y, t) := 
\sigma_{(x(t),y(t)), h(t)}(x, y, t),~~~~~\mbox{with}~~h(t) = bt,
$$
and is the identity map on the
boundary of the cone.

%
%

\begin{lemma}~\label{tube.curve}
Suppose  $0 < b \le 1$.
The map $\tau: C(b, h)  \to C(b, h)$ 
associated to the curve $\alpha(t)$ is a homeomorphism of $C(b,h)$,
leaving its boundary fixed. 
This map is a Lipschitz map with Lipschitz constant $C_{L}(\tau)$
satisfying
$$
C_{L}(\tau) \le  4\sqrt{1 + b^2}  + C_L(\alpha) + 2.
$$
\end{lemma}
{\bf Proof:}
Let $u= (x_0, y_0, t_0), v= (x_1, y_1, t_1)$ be two points in $C(b,h)$.
If they lie at the same height $t_0=t_1$,
then since $\alpha(t_0)$ lies inside  $C_L$, 
Lemma~\ref{sigma.lipschitz} yields the Lipschitz bound
$ ||\tau(u) - \tau(v)|| \le 2||u - v||$.

Now if $t_0 \ne t_1$
we may arrange that $t^{\ast} :=t_1 - t_0 > 0$,
by exchanging $u$ and $v$ if necessary. Set $B := \sqrt{1 + b^2},$
and our  object is to show that 
$$
||\tau(u) - \tau(v)|| \le (4B  + C_L(\alpha) + 2) ||u - v||.
$$

We study the effect of radial projection using  the cone
vertex $(0,0,0)$. Given a point $w=(x, y, t_0)$ in the
plane $z=t_0$ we let $P(w) = (t_1/t_0) w$ be its radial projection
to the plane $t= t_1$. Then for all such points $w$ in the cone
$C(g, h)$  we have
$$
t^{\ast} \le ||P(w) - w|| \le (\sqrt{1+ b^2}) t^{\ast}= B t^{\ast},
$$
the upper bound being tight for $w$ on  the edge of the cone.

Since $P(u)$ and $v$ both lie in the cone in the
same plane $t=t_1$, we have 
\begin{eqnarray*}
||\tau(u) - \tau(v)|| & \le &  || \tau(u) - P(\tau(u)) || +
||P(\tau(u))- \tau(P(u))|| +
||\tau(P(u)) - \tau(v)|| \\
& \le & B t^{\ast} + ||P(\tau(u))- \tau(P(u))||
+ 2||P(u) - v||.
\end{eqnarray*}
To bound the  term $||P(\tau(u))- \tau(P(u))||$
we observe that in terms of the
function of Lemma~\ref{sigma.lipschitz} applied in the
plane $t= t_1$ of the cone $C(b,h)$ we have 
$P(\tau(u)) = \sigma_{ P(\alpha(t_0)), bt_1}(P(u))$
while $\tau(P(u))= \sigma_{ \alpha(t_1), bt_1}(P(u))$. Then
\begin{eqnarray*}
||\sigma_{ P(\alpha(t_0)), bt_1}(P(u)) - 
\sigma_{ \alpha(t_1), bt_1}(P(u))||
& \le & ||P(\alpha(t_0)) - \alpha(t_1)|| \\
& \le & ||  P(\alpha(t_0))  - \alpha(t_0)|| + ||\alpha(t_0) -  \alpha(t_1)|| \\
& \le & B t^{\ast}  + C_L(\alpha)t^{\ast},
\end{eqnarray*}
where the first line above follows using
Lemma~\ref{sigma.lipschitz}.
Substituting this in the previous inequality yields
$$
||\tau(u) - \tau(v)|| \le \left(2B + C_L(\alpha)\right) t^{\ast} 
+ 2||P(u) - v||.
$$
On the other hand we have
$||u - v|| \ge t$ and 
$$
||u -v || \ge ||P(u) - v|| - ||u - P(u)||
\ge ||P(u) - v|| - Bt^{\ast}.
$$
We average these two inequalities in a suitable  ratio to obtain
\begin{eqnarray*}
||u - v|| &\ge & \frac{4B + C_L(\alpha)}
{4B + C_L(\alpha) +2} t^{\ast}  + \frac{2}
{4B + C_L(\alpha) +2}(||P(u) - v|| - Bt^{\ast} ) \\
 & \ge & \frac{1}{4B + C_L(\alpha) +2}
\left( (2B + C_L(\alpha)) t^{\ast}  
+ 2||P(u) - v||\right) \\
& \ge &  \frac{1}{4B + C_L(\alpha) +2}||\tau(u) - \tau(v)||,
\end{eqnarray*}
which gives the Lipschitz bound for $\tau$. $~~~\qed$ \\

%
%

{\bf Proof of Theorem~\ref{unknot.smoothupper}:}
The curve $\gamma$ has a normal tubular neighborhood
of radius $r.$ As described above 
we inscribe in $\gamma$  a closed polygon $P'$ with $n$ edges,
with $\ds n \le \frac{32L}{r}$, labelling the
points $z_j = \gamma(s_j)$, the  
$j$-th edge $e_j= [z_j, z_{j+1}]$ of $P'$,
and   $\gamma_j$  the corresponding
arc of the curve between $\gamma(s_j)$ and 
$\gamma(s_{j+1})$, as described above. We also have
the set of $n$ right circular cones $\{C_j: 0 \le j \le n-1\}$, where
$C_j$ has the segment $e_j$ as its axis.
Now Lemma~\ref{cones.disjoint} applies to show
that $P'$ is embedded
in $\RR^3$ and 
that these cones are disjoint except at the points $z_j$.

We define a homeomorphism $\tau: \RR^3 \to \RR^3$,
which is the identity outside all the cones $C_j$ and
which inside each cone $C_j$ is given by the homeomorphism 
$\tau_j$ of $C_j$
given by Lemma~\ref{tube.curve}, 
taking the  arc $\alpha(t)$
to be the arc $\gamma_j$, reparameterized using 
the variable $t$ that
measures distance along $e_j$, rather than
the arc length variable $s$.  With this reparameterization
the map $t \to \alpha(t)$ is a Lipschitz map.
To bound the Lipschitz constant we use 
results from the proof of Lemma~\ref{cones.disjoint}.
It shows that the tangent vector $T$ in the
direction of $e_j$ satisfies  $1 \ge T \cdot T(s')  \ge  \cos \theta_a$
for all tangent vectors on the arc $\gamma_j$, 
where $0< \theta_a \le 8 \arcsin 1/64$ 
which yields $\cos \theta_a > 0.9$. This implies
that the Lipschitz constant $C_L(\alpha) \le 2.$

We conclude using Lemma~\ref{tube.curve} that
$\tau$ is a Lipschitz map and obtain a bound
for its Lipschitz constant inside the individual
cones $C_j$. We take 
$h = \mbox{length}(e)$ and 
$b \le \sin (8 \arcsin (\theta_{a_j}/2)) \le \sin (0.14) < 0.14,$
which yields the bound
$C_L(\tau) \le 4(0.14) + 2 + 2\le 5.$ Since $\tau$ is
the identity map outside the cones, this Lipschitz
bound carries over to $\tau: \RR^3 \to \RR^3$, using
the triangle inequality.

Now consider the triangulated PL surface $\Sigma'$ having
$P'$ as boundary given by Theorem~\ref{PLupperbound},
having  area $A'$ bounded by
$$
A' < (C_4)^{n^2}(L')^2,
$$
where $L'$ is the length of $P'.$ Using the bound $n \le 33L/r$
and that $L' \leq L$ and  $\ds L/r >1,$ this gives
$$
A' < (C_4)^{(33L/r)^2}(L')^2 \le (C_0)^{(L/r)^2} L^2,
$$
for an enlarged constant $C_0 = (C_4)^{33}$.

We can assume, possibly after a small perturbation
 that induces an arbitrarily
small increase in the area of $\Sigma'$, that $\Sigma'$
intersects the boundary of the $r/2$-tubular neighborhood of $\gamma$
and the boundary of each ball $B_j$ transversely.
This follows from the fact that
the space of smooth disks is dense in the
space of piecewise smooth disks in $\RR^3$, and
the area functional is continuous on the space of piecewise smooth disks 
with the smooth topology.
See for example Kosinski \cite[Theorem 2.5]{K93}. \\

Applying the homeomorphism $\tau$ then gives a piecewise-$C^2$
embedded surface $\Sigma := \tau(\Sigma')$ that has $\gamma$
as boundary and whose area is increased
from the area $A$ of $\Sigma'$ by at most a factor of  25,
using the bound of 5 on the Lipschitz constant of the map $\tau$.
(A Lipschitz map
$\tau: \RR^n \to \RR^n$ with Lipschitz constant
$C_L(\tau)$ 
can stretch the volume of an $m$-dimensional submanifold
by at most a factor $C_L(\tau)^m$, see
Federer\cite[2.10.11]{Fe69}, Morgan~\cite[Prop. 3.5]{Mo88}.)
The fact that $\Sigma$ is piecewise
$C^2$ follows using
the fact that the map $\tau$ is piecewise $C^2$, with
its smoothness inherited from the smoothness of the curve $\gamma(s)$.
Thus, enlarging the constant
$C_0$ to a suitable $C_6$ we have
$$
A < \frac{1}{2}C_6^{(L/r)^2} L^2.
$$

Finally, we can perturb the interior of $\Sigma$ to obtain a smooth
surface while again changing its area by an arbitrarily small amount,
and without moving the boundary $\gamma$, see \cite[Theorem 2.5]{K93}. 
This yields the theorem, where we give up another factor of $2$ on
the right side of the bound above to absorb the effect of the perturbation.
$~~~\qed$

%
%
%
%

\section{Lower bound: Large complexity minimal
PL spanning disks}\label{PL}

In this section
we construct a series of polygonal unknotted curves $K_n$
in $\RR^3$ whose minimal spanning disks must have exponential
complexity, according to various measures. This is an extension
of the results of \cite{HST}.

%
%

\begin{theorem} \label{properties}
There exists a sequence of unknotted  polygonal curves $\{K_n\}$
in $\RR^3$, for $n \ge 1$,  having the following properties.

\noindent
(1) $K_n$ contains at most $24n+40$ edges.

\noindent
(2) Any embedded, piecewise-smooth disk
spanning $K_n$ intersects the $y$-axis
in at least $2^{n-1}$ points.

\noindent
(3) Any embedded PL triangulated disk $D_n$ bounded by $K_n$
contains at least $2^{n-1}$ triangles.

\noindent
(4) Each $K_n$ has length one and
there exists a constant  $r_0 > 0$ independent of $n$ such that each
$K_n$ is disjoint from the solid
cylinder $Cyl(r_0)$ of radius $r_0$ around the $y$-axis.

\end{theorem}

\begin{figure}[hbtp]
\centering
\begin{minipage}[c]{.3 \textwidth}
\centering
\includegraphics[width=.3\textwidth]{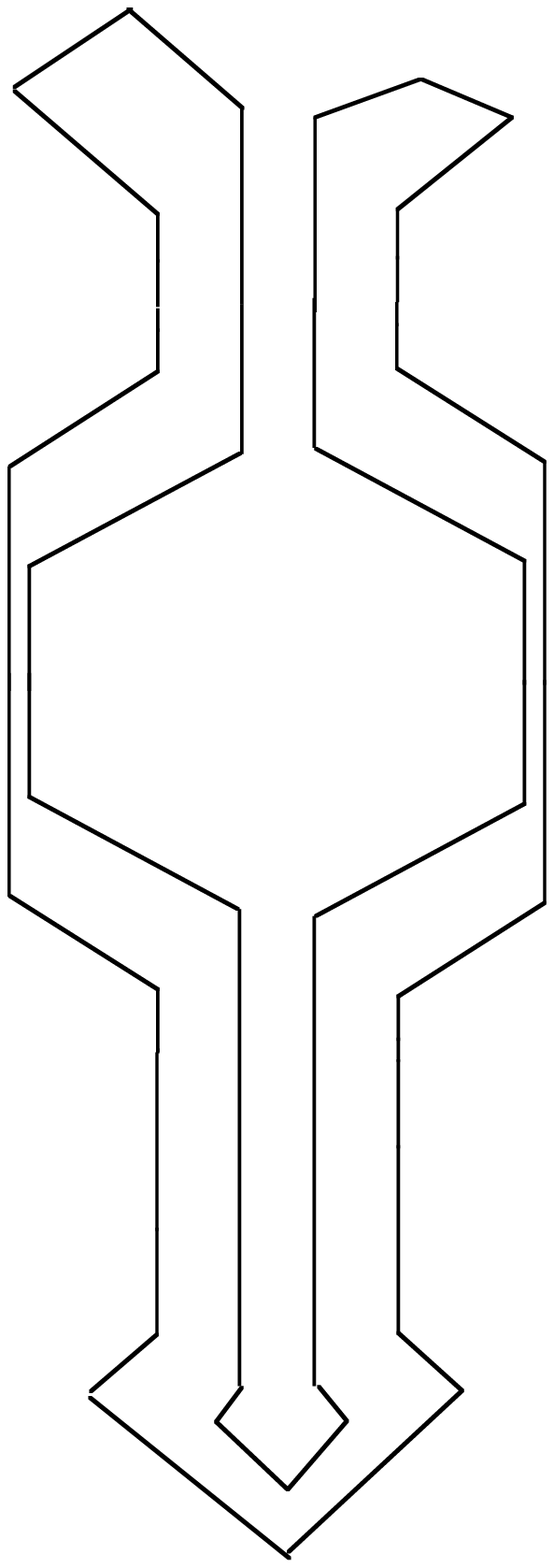}
\end{minipage}
\begin{minipage}[c]{.3 \textwidth}
\centering
\includegraphics[width=.26\textwidth]{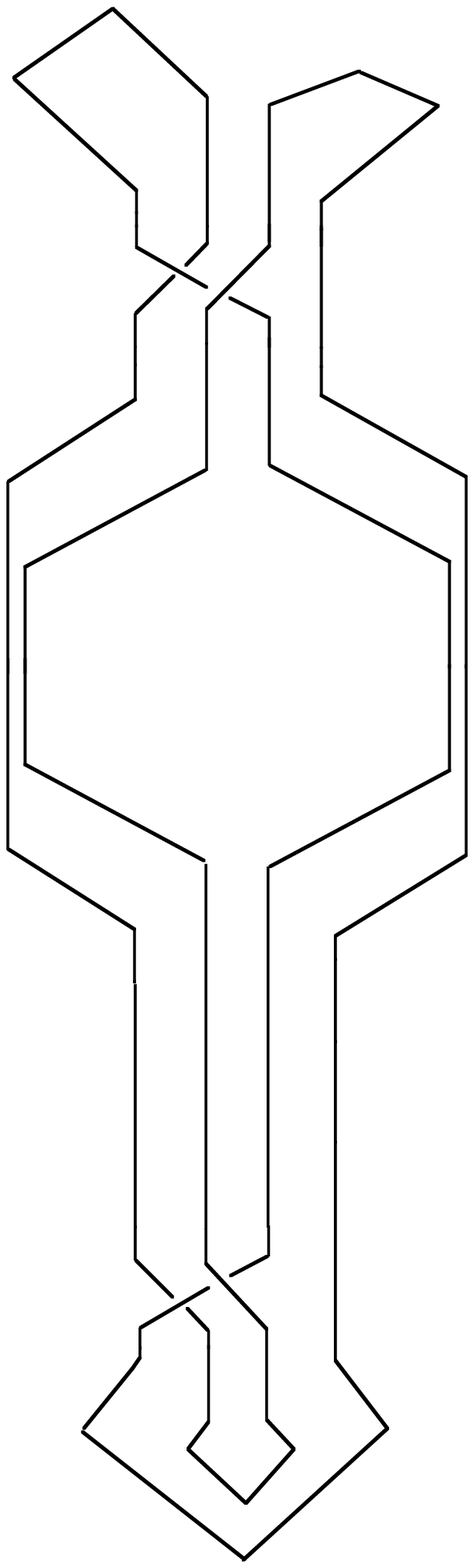}
\end{minipage}
\begin{minipage}[c]{.3 \textwidth}
\centering
\includegraphics[width=.24\textwidth]{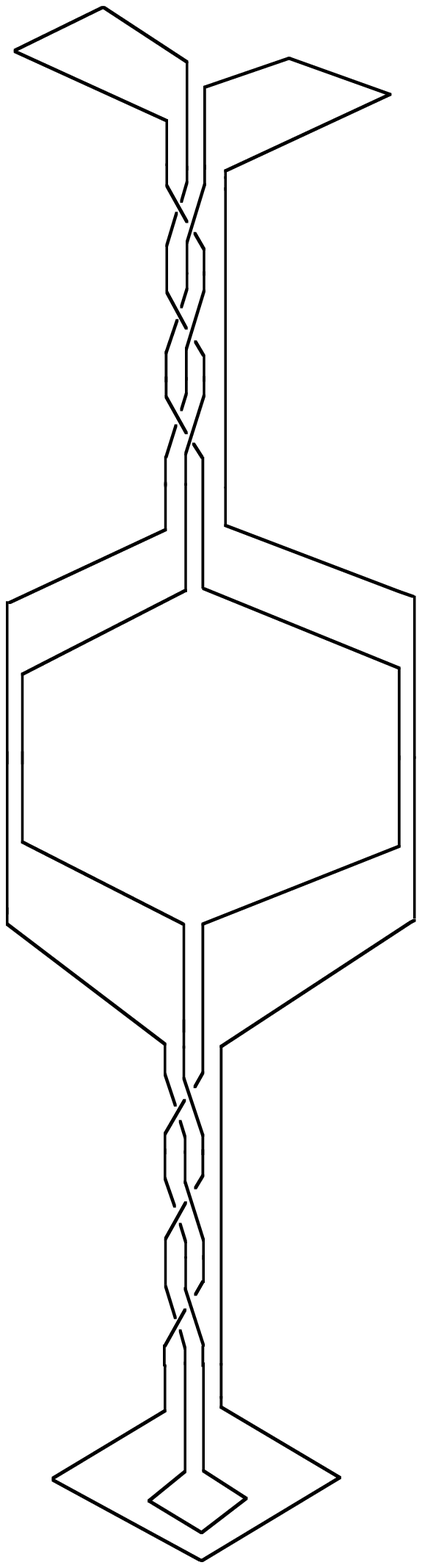}
\end{minipage}
\caption{The curves $K_0$, $K_1$ and $K_3$.}
\label{Kn}
\end{figure}

Three curves in this sequence are shown in Figure~\ref{Kn}.

%
%

\noindent\paragraph{Remark.} The construction of the curves $K_n$
in Theorem~\ref{properties}
differs from that in \cite{HST}, in order to make
condition (4)  hold.
The main focus of \cite{HST} was to construct examples of curves
with $O(n)$ edges satisfying property (3), while in the
application here it is properties (2) and (4) that are needed.

%
%

\begin{figure}[hbtp]
\centering
\includegraphics[width=.4\textwidth]{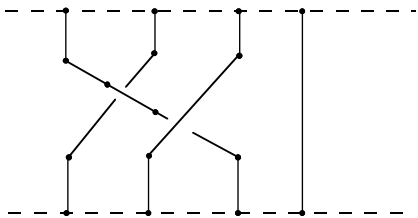}
\caption{The braid $\alpha = \sigma_1 \sigma_2 ^{-1}$ can be embedded in
$\RR^3$ using 12 segments so that the top and bottom ends of each arc
are vertical segments on the $xz$-plane. The leftmost strand is made from
five segments,
and rises above and below the  $xz$-plane.
The other three strands lie completely in the  $xz$-plane.}
\label{alpha}
\end{figure}

\noindent
{\bf Proof:}
The assertions (1)-(3) of the theorem are invariant under homothety
$\bx \to \lambda \bx$,
so restricting to curves $K_n$ of length one results in
no loss of generality.

\begin{figure}[hbtp]
\centering
\includegraphics[width=.15\textwidth]{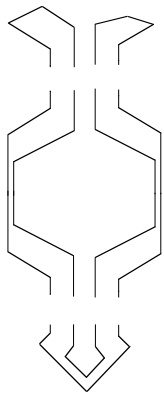}
\caption{Arcs $a_n, b_n, \dots, h_n$ used in forming the top, middle and
bottom pieces of $K_n$.
As $n$ increases, the horizontal separation between the endpoints of the
arcs
decreases to 0, but their total length remains equal to 1/3.
All arcs lie in the $xz$-plane
and are disjoint from a disk in this plane
of radius $r_0$ centered at the origin.}
\label{caps}
\end{figure}

The construction of $K_n$
begins with the 4-braid $\alpha = \sigma_1 \sigma_2 ^{-1}$
shown in Figure~\ref{alpha}.
We construct $\alpha$ to consist of vertical segments near
its endpoints and to have (total) length one. We
stack $n$ copies of $\alpha$ to get a braid $\alpha^n$, and then rescale
by a homothety to get a braid
$\displaystyle \alpha_n = \frac{1}{3n}\alpha^n$ whose length is $1/3$.
Similarly we stack $n$ copies of $\alpha^{-1}$ and scale
to get a braid $\alpha_{-n}$ of length $1/3$.
The curve $K_n$ is constructed by starting with a copy of
each of $\alpha_n,\  \alpha_{-n}$ and adding eight arcs to obtain a
closed curve in $\RR^3$. The eight arcs are shown in Figure~\ref{caps}.

%
%

\begin{figure}[hbtp]
\centering
\includegraphics[width=.25\textwidth]{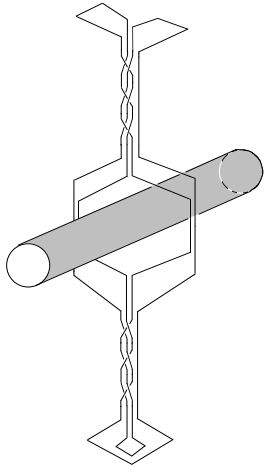}
\caption{Each $K_n$ is disjoint from a cylinder of
radius $r_0$ around the $y$-axis.}
\label{cylinder}
\end{figure}

The added eight arcs all lie in the $xz$-plane, as shown
in Figure~\ref{caps}. They consist of two arcs $a_n, b_n$
connecting the first and second and the
third and fourth endpoints at the top of $\alpha_n$, four arcs
$c_n,d_n,e_n,f_n$ connecting the bottom four endpoints of  $\alpha_n$
to the top four
endpoints of $\alpha_{-n}$, and two arcs $g_n, h_n$ connecting the
first and fourth and the second and third endpoints of the
bottom of  $\alpha_{-n}$.
Forty line segments are used to construct  $a_n, b_n, c_n,d_n,e_n,f_n,
g_n, h_n$.
The resulting closed curves $K_n$, shown in Figure~\ref{Kn},
are each unknotted, since the braid $\alpha^{n}\alpha^{-n}$ is trivial, and
each $K_n$ uses at most $24n + 40$ segments, with a contribution of
$12n$ coming from each of $\alpha_{n}$, $\alpha_{-n}$.
For each $n$, we construct the eight added arcs so that they have total
length $1/3$, and so that they are disjoint from
a solid cylinder of fixed radius $r_0$.
With the total length of all the segments set to
one, a calculation shows that
we can take $r_0$ to be $1/100$. See Figure~\ref{cylinder}.

%
%

As $n \to \infty$ the separation between the top four points
of $\alpha_{n}$ converges to zero, and the
endpoints of $a_n$ and $b_n$ which are identified to them also converge.
Endpoints of $c_n,\dots , h_n$ converge similarly.
Each curve $K_n$ is a length one curve disjoint from the
cylinder $Cyl(r_0)$ of radius $r_0$ around the $y$-axis.

Now consider an arbitrary piecewise smooth disk $D_n$ with
boundary $K_n$. The intersections of $D_n$ with the level sets
$\{ z = c \}$ can be analyzed using Morse Theory and the theory
of train tracks and surface diffeomorphisms.
This analysis is carried out in Theorem~1 of \cite{HST},
where it is shown that the $y$-axis must intersect
$D_n$ in at least  $2^{n-1}$ points.
All cases of Theorem~\ref{properties} have now been established.
$\qed$

%
%
%

\section{Lower Bound: Large minimal area
embedded spanning disks}~\label{lower.bounds}

In this section we show that
there is no upper bound
on the minimal area of any embedded disk that depends
only on the length $L$ of the curve,
as stated in Theorem~\ref{noupperbound2}.
This is based on
the sequence of curves $K_n$ constructed in Theorem~\ref{properties}.

%
%
 
\begin{lemma} \label{perturb}
Let $\{K_n\}_{n=1,...,\infty}$ be the sequence of length one curves
constructed in  Theorem~\ref{properties}.
There is a constant $c_7 > 1$ such that any
embedded disk $D_n$ having boundary $K_n$ has area $A$ satisfying
$$A > c_7 2 ^{n}.$$
The same lower bound holds for 
the area of any embedded disk having boundary
a curve $\gamma'$ that is isotopic to $K_n$ by an isotopy
that fixes the cylinder $Cyl(r_0)$ of radius $r_0$ around the $y$-axis.
\end{lemma}

{\bf Proof:}
Suppose a line $L$ in $Cyl(r_0)$ parallel to its axis intersects $D_n$ in
less than $2^{n-1}$ points. There is an isotopy of $\RR^3$ fixing
$K_n$ and carrying $L$ to the $y$-axis.  The image of $D_n$ under this
isotopy
would be a disk $D_n'$ spanning $K_n$ and intersecting the $y$-axis
in less than $2^{n-1}$ points, contradicting Theorem~1 of \cite{HST}.
So any such line in the cylinder $Cyl(r_0)$ intersects $D_n$ in
at least $2^{n-1}$ points.
The area of a perpendicular cross-section of $Cyl(r_0)$ is
$\pi r_0 ^2 $, so the projection of $D_n$ to the $xz$-plane has
area at least $2^{n-1} \pi r_0 ^2 $, as does $D_n$ itself.
We choose $c_7=\pi r_0 ^2 /2.$

Now suppose that $L_n$ is a curve isotopic to $K_n$ in the complement of
$Cyl(r_0)$ and $E_n$ is a disk spanning $L_n$.  Then there is an isotopy of
$\RR^3$ carrying $E_n$ to a disk spanning $K_n$ that restricts to the
identity on $Cyl(r_0)$. The area of $E_n$ in  $Cyl(r_0)$ is at least
$c_7 2 ^{n}$  by the previous argument.
$ \qed $

\begin{figure}[hbtp]
\centering
\begin{minipage}[c]{.49 \textwidth}
\centering
\includegraphics[width=.25\textwidth]{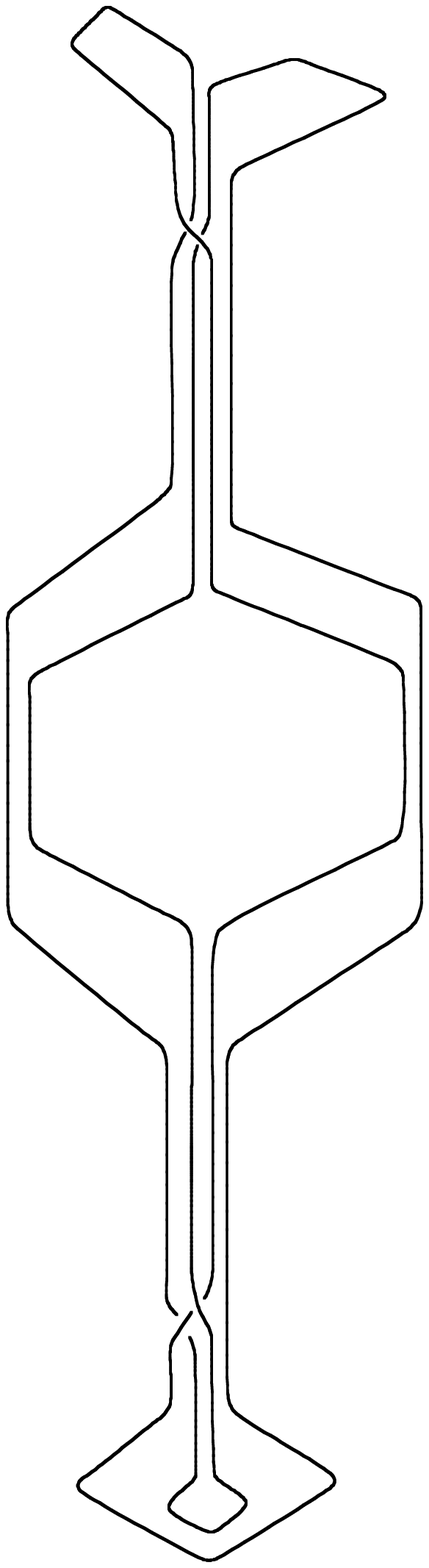}
\end{minipage}
\begin{minipage}[c]{.49 \textwidth}
\centering
\includegraphics[width=.25\textwidth]{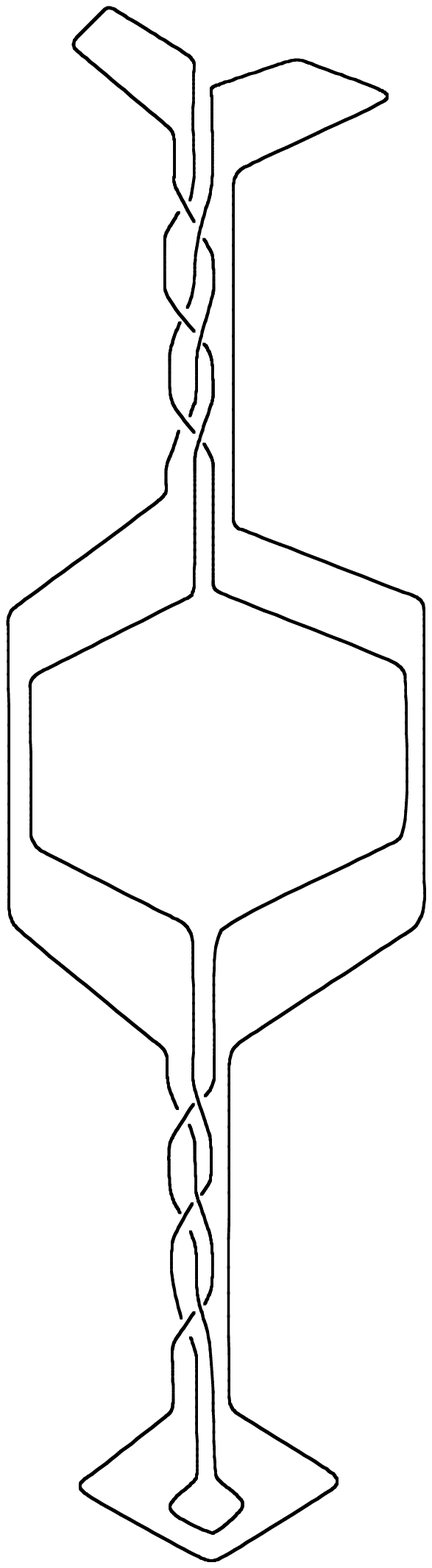}
\end{minipage}
\caption{Each $\tilde{\Delta}_n$ is a smooth curve that
has length one and thickness at least $k/n$. $\tilde{\Delta}_1$
and $\tilde{\Delta}_3$ are shown.}
\label{Deltan}
\end{figure}

To obtain $\Delta_n$ we will construct
a smooth curve of length one that approximates $K_n$,
and in particular remains disjoint from $Cyl(r_0)$ and has
thickness greater than $1/n$,
as shown in Figure~\ref{Deltan}. \\

%
%

{\bf Proof of Theorem~\ref{noupperbound1}:}
It suffices to find curves $\Delta_n$ having
length at most one with the other claimed properties, as rescaling
by a homothety to make these curves have length
precisely one increases both the
thickness and the area of any spanning disk.
We will approximate each $K_n$ by a smooth
curve $\tilde{\Delta}_n$,
isotopic to $K_n$ in the complement of the cylinder $Cyl(r_0)$.
Finally we will take $\Delta_n := \tilde{\Delta}_{[cn]+1}$ for an
appropriate constant $c\ge 1$.

We first smooth the braids $\alpha$ and $\alpha^{-1}$ by replacing
neighborhoods of their vertices by smooth arcs, getting a smooth
braid $\beta$ of thickness $k > 0$ and
curvature bounded above by a constant $K_0$. We do this smoothing
so that $\beta$ still consists of vertical arcs near its endpoints.
Then there is a constant $r_1 > 0$ with the property that the
thickness of $\beta$ is greater than $r_1$, and each vertical segment
near an endpoint of $\beta$ has length greater than $r_1$.
Form $\beta^n$ and $\beta^{-n}$
as before by stacking copies of $\beta$ or $\beta^{-1}$. Since copies
of $\beta$ are identified along points that have neighborhoods coinciding
with vertical segments, the resulting curve is $C^\infty$.
Moreover the thickness of $\beta^n$ and $\beta^{-n}$ is still greater
than $r_1$, since distinct copies of
$\beta$ have $r_1$-tubular neighborhoods with disjoint interiors.
Scaling by a homothety gives
$\displaystyle \beta_n = \frac{1}{3n}\beta^n$, with length $1/3$, and
similarly $\beta_{-n}$.
The thickness of each of $\beta_n$ and  $\beta_{-n}$ is greater than
$\displaystyle \frac{r_1}{3n}$.
We can also smooth the curves $a_n, \dots , h_n$ in neighborhoods
of their internal vertices,
to obtain nearby arcs  $a_n', \dots , h_n'$ that are smooth
with arcs near the endpoints remaining as straight vertical segments,
with the length of each straight segment no less than a constant $r_2$
(which does not depend on $n$).
The approximating arcs are chosen sufficiently close
to maintain their disjointness from  $Cyl(r_0)$, and with their
total length remaining less than $1/3$.
Then combining $ \beta_n, \beta_{-n}$ with
$a_n', \dots , h_n'$ gives a closed curve
$\tilde{\Delta}_n$ of length one that is
isotopic to $K_n$ in the complement of $Cyl(r_0)$, and the
thickness of $\tilde{\Delta}_n$ is given by
$\displaystyle \frac{k}{3n}$, where $k$ is the smaller of $r_1, r_2$.

%
%

We now establish a lower bound on the area $A$ of any
embedded spanning disk for $\tilde{\Delta}_n$.
Lemma~\ref{perturb} implies that
there is a $C_0 > 1$ such that for large enough $n$, any
embedded disk with boundary $\tilde{\Delta}_n$ has area $A$
greater than ${C_0}^n$.

Note that $k<1$ since $\tilde{\Delta}_n$ has length one
and thickness less than
$k$. Thus we can
rescale $n$ to eliminate the constant $k$
appearing in the thickness bound
by taking a subsequence $\Delta_n :=\tilde{\Delta}_{[kn] +1}$,
and setting  $C_1 = (C_0)^k.$
We then have minimal area $A > (C_1)^n$.
$\qed$ \\

In the next result
   we show  how to deform $K_n$ to construct a similar sequence
of smooth curves $\Gamma_n$ having uniformly bounded curvature.
Recall that the curvature function for a $C^2$-curve
embedded in $\RR^3$ is computed by taking
an arclength parameterization $\gamma: [0, L] \to \RR^3$,
so that $\displaystyle ||\frac{d}{dt} \gamma(t)||=1$, and the curvature function is
then $\displaystyle \kappa(t) = ||\frac{d^2}{dt^2}\gamma(t)||.$ \\

%
%

{\bf Proof of Theorem~\ref{noupperbound2}:}
The construction of the family
$\{\Gamma_n\}_{n=1,...,\infty}$ is similar to that of
 $\{\Delta_n\}_{n=1,..., \infty}$.
We again use smoothings $a_n', \dots , h_n'$ of
the arcs $a_n, \dots , h_n$ to curves
that lie nearby them, and assert that the approximating curves
can be chosen to have
curvature at each point bounded above by a fixed constant $K$.
To see this we use the fact that a sequence of smooth curves
converging to a smooth limit curve have curvatures converging
pointwise to the curvature of the limit curve.
For any curve with monotonically increasing
$z$-component, a deformation of the curve by a diffeomorphism
$g_{\lambda}(x,y,z) := (\lambda x, \lambda y, z)$
carries it to a curve that, for $\lambda$ small, lies very close to
to the $z$-axis.
For $\lambda$ sufficiently small the
curvature of the image curve under $g_{\lambda}$
is uniformly close to zero.
We choose $\lambda_n \to 0$ to be a sequence converging to 0 sufficiently
fast so that  the curvature of
$\beta_n  = g_{\lambda_n}(\beta^n) \to 0$, where $\beta^n$
is as before, and similarly we define $\beta_{-n}$.
We also smooth the curves $a_n, \dots , h_n$ in neighborhoods
of their internal vertices,
to obtain nearby arcs  $a_n', \dots , h_n'$ that are smooth
with arcs near the endpoints remaining straight vertical segments.
While these vertical segments all converge to the $z$-axis as $n \to
\infty$,
the curvature of each of $a_n', \dots , h_n'$ remains uniformly bounded.
One way to see this explicitly is to choose smooth limiting curves
$a_\infty', \dots , h_\infty'$
of bounded curvature that are not embedded, but rather have arcs near
their endpoints that agree
with the $z$-axis. Any sequence of smooth curves $a_n', \dots , h_n'$ that
converge to these curves smoothly
will have uniformly  bounded curvature.
Then the curvature of the curves
$\Gamma_n$ actually converges to 0 at points corresponding
to the braids  $\beta^n$, $\beta_{-n}$, and is uniformly bounded at other
points.
The curvature is zero along the straight segments where
$a_n', \dots , h_n'$ join $\beta^n$ and $\beta_{-n}$.

Since the required area for a spanning disk goes to infinity as
$n \to \infty$, we can rechoose $\Gamma_n$ to be a suitable
subsequence to have the required property $A > n.$
This concludes the proof of Theorem~\ref{noupperbound2}.
$~~~~~ \qed$

%
%

\section{Appendix: Isoperimetric theorems for curves in
$\RR^3$}~\label{classical}

This appendix establishes 
two versions of the isoperimetric
theorem in $\RR^3$ as stated as
(1) and (2) in \S1. Result (1)
follows from the solution of the
Plateau problem, due to Douglas and Rado, \cite{Douglas}, \cite {Rado}, and
on a result of Carleman on isoperimetric inequalities for
minimal disks \cite{Car21}.  See also Reid \cite{Reid}.
We  deduce (2) from (1),
using  standard cut and paste methods of 3-manifold topology,
which can be used to convert
an immersed surface to an embedded surface, possibly
with an increase in genus, but maintaining orientability.
Such a cut and paste can be achieved without
increasing area by a ``rounding of creases''.
Alternate approaches to proving (2)
are possible based on work
of Blaschke \cite[p. 247]{Bla30},
which was the original approach to (2).
The discussion in Osserman \cite[p. 1202]{Os78} indicates that Blaschke's
argument is heuristic, but
can be made valid for area-minimizing immersed surfaces
having the curve as boundary. \\

\begin{theorem}~\label{iso-immersed}
Let $\gamma$ be a simple closed $C^2$-curve in $\RR^3$
of length $L$. Then there exists an immersed disk
having area $A$ in $\RR^3$, with $\gamma$ as boundary,
such that
$$
4 \pi A \le L^2.
$$
If the curve is not a circle, then there exists
such an immersed disk for which strict inequality holds.
\end{theorem}

{\bf Proof of  Theorem~\ref{iso-immersed}:}
Span $\gamma$ by a least area disk $D$, whose
existence is guaranteed for rectifiable $\gamma$ by the solution
of the Plateau problem \cite{Douglas}, \cite{Rado}.
The regularity results of Osserman \cite{Os70} and
Gulliver \cite{Gulliver} show that $D$ is a smooth immersion
in its interior. (The interior of $D$ may transversely
intersect $\gamma$, as indeed it must if $\gamma$ is knotted.)
The argument in Reid \cite{Reid} (see also \cite{Car21}) shows that
$4 \pi A \le L^2$, with equality if and only if $\gamma$ is a circle.
$~~~\qed$ \\

\begin{theorem}~\label{iso-embedded}
Let $\gamma$ be an embedded closed $C^2$-curve in $\RR^3$
of length $L$. Then there exists an embedded orientable smooth
spanning surface for $\gamma$ in $\RR^3$ having area $A$, with
\begin{equation}
   4 \pi A \le L^2.
\label{isopineq}
\end{equation}
If the curve is not a circle, then there exists
such a surface for which strict inequality holds.
\end{theorem}

{\bf Proof of  Theorem~\ref{iso-embedded}:}
Apply Theorem~\ref{iso-immersed} to get an immersed disk $D$ spanning
$\gamma$ satisfying the isoperimetric inequality (\ref{isop}).
If $D$ is not embedded then
perturb it slightly
so that its self-intersections are in general position, i.e. transverse
with finitely many triple points.
Since smooth general position maps are dense among $C^2$ maps, this can be
done with an arbitrarily small increase in area.
The self-intersections then consist of a finite number of
immersed double curves.  There is
a uniquely defined surface $F$ obtained by oriented cut-and-paste along
these curves.
An arbitrarily small perturbation of $F$
moves it to an embedded surface.
The resulting surface is orientable
with the same boundary as $D$,
but may have higher genus. If $D$ has area strictly less
than $L^2 / 4 \pi$ then the area of $F$ can also be taken to be less than
$L^2 / 4 \pi$.
If $D$ has area equal to $L^2 / 4 \pi$ then it spans a circle
and is already embedded. In all other cases strict inequality applies
for both $D$ and $F$.
$~~~\qed$ \\

\noindent
{\bf Remark:}
Almgren  \cite{Alm86} established
isoperimetric inequalities for varifolds and currents
in arbitrary codimension, with optimal constants.
In some range of dimension and codimension, these give upper
bounds for the area of embedded spanning surfaces and spanning
submanifolds, but with unspecified topological type.

\vspace{+10 pt}
\noindent
Department of Mathematics \\
University of California \\
Davis, CA 95616 \\
hass@math.ucdavis.edu \\

\vspace{+10 pt}
\noindent
AT\&T Labs-Research \\
Florham Park, NJ 07932-0971 \\
jcl@research.att.com \\

\vspace{+10 pt}
\noindent
Department of Mathematics \\
University of California \\
Davis, CA 95616 \\
wpt@math.ucdavis.edu
\end{document}